\documentclass[12pt]{article}
\usepackage{amssymb}
\usepackage[dvips]{graphicx}
\usepackage{latexsym}
\usepackage{amsmath}
\usepackage{amssymb,color}
\usepackage{amsthm}
\usepackage{url}
\usepackage{xspace}
\usepackage{paralist}

\newtheorem*{thm*}{Theorem}

\theoremstyle{definition}

\theoremstyle{remark}

\input{tcilatex}
\begin{document}

\title{A New Topological Helly Theorem and some Transversal Results}
\author{L.~Montejano \thanks{%
The author acknowledges support of CONACYT, proyect 166306} \\
%EndAName
{\small Instituto de Matem\'{a}ticas, Unidad Juriquilla. }\\
{\small Universidad Nacional Aut\'{o}noma de M\'{e}xico, M\'{e}xico} }
\maketitle

\begin{abstract}
\noindent We prove that for a topological space $X$ with the property that $%
H_{\ast}(U)=0$ for $\ast \geq d$ and every open subset $U$ of $X$, a finite
family of open sets in $X$ has nonempty intersection if for any subfamily of
size $j,$ $1\leq j\leq d+1,$ the $(d-j)$-dimensional homology group of its
intersection is zero.

\noindent We use this theorem to prove new results concerning transversal
affine planes to families of convex sets \ \ 
\end{abstract}

%%%%%%%%%%%%%%%%%%%%



\section{\noindent Introduction and Preliminaries}

\noindent A prominent role in combinatorial geometry is played by the Helly
Theorem \cite{DGK}, which states that a finite family of convex sets in $%
\mathbb{R}^{d}$ has nonempty intersection if and only if any subfamily of
size at most $d+1$ has nonempty intersection. Helly himself realized in 1930
(see \cite{E}) that a finite family of sets in $\mathbb{R}^{d}$ has nonempty
intersection if for any subfamily of size at most $d+1,$ its intersection is
homeomorphic to a ball in $\mathbb{R}^{d}$. In fact, the result is true if
we replace the notion of topological ball by the notion of acyclic set, see 
\cite{Bog} and \cite{K}. \ \ In 1970, Debrunner \cite{Deb} proved that a
finite family of open sets in $\mathbb{R}^{d}$ has nonempty intersection if
for any subfamily of size $j,$ $1\leq j\leq d+1,$ its intersection is $(d-j)$%
-acyclic. In fact, these hypothesis imply that each one the open sets are
acyclic.

Unlike all previous topological Helly results in the literature, where the
ambient space is the euclidean space and the sets are acyclic, we only
require as an ambient space a topological space $X$ in which $H_{\ast }(U)=0$
for $\ast \geq d$ \ and for any subfamily of size $j$ of our sets, the $%
(d-j) $-dimensional homology group of its intersection to be zero.

The fact that this is a non-expensive topological Helly theorem \ ---in the
sense that it does not require the open sets to be simple--- from the
homotopy point of view (we only require its $(d-1)$-dimensional homology
group to be zero), allows us to prove interesting new results concerning
transversal planes to families of convex sets.

During this paper, we use reduced singular homology with any nonzero
coefficient group. \ Let $U$ be a topological space. We say that $%
H_{-1}(U)=0 $ if and only if $U$ is nonempty and we say that $U$ is \textit{%
connected } if and only if $H_{0}(U)=0$. $\ $\ For $q=-1$ and $q=0$, the
exactness of the Mayer-Vietoris sequence should be understood as the
statement that if $\ H_{q}(U)$ occurs there between two vanishing terms,
then $H_{q}(U)=0.$ \ For an integer $n\geq -1,$ we say that $U$ is $n$-%
\textit{acyclic} if $H_{\ast }(U)=0$ for $-1\leq \ast \leq n.$ Furthermore,
\ $U$ is \textit{acyclic} if $H_{\ast }(U)=0$ for $\ast \geq -1.$

\bigskip

\section{A Topological Helly-Type Theorem}

\smallskip We start with the following auxiliary proposition.

\smallskip

\noindent \textbf{Proposition A}$_{(m,\lambda )}.$\ \textit{Let }$%
F=\{A_{1},...,A_{m}\}$\textit{\ be a family of open subsets of a topological
space }$X$\textit{\ and let }$\lambda \geq 0$\textit{\ be an integer.
Suppose that for any subfamily }$F^{^{\prime }}\subset F$\textit{\ of size }$%
j,$\textit{\ }$\ 1\leq j\leq m$\textit{,}

\begin{center}
$H_{m-1-j+\lambda }(\bigcap F^{^{\prime }})=0$.
\end{center}

\textit{Then}

\begin{center}
$\mathit{\bigskip }H_{m-2+\lambda }(\bigcup F)=0$.
\end{center}

\noindent \textbf{Proof. \ }Proposition A$_{(2,0)}$ claims that the union of
two connected sets with nonempty intersection is a connected set, and
Proposition A$_{(2,\lambda )}$ is just the statement of the exactness of the
Mayer-Vietoris sequence: \ $0=H_{\lambda }(A_{1})\oplus H_{\lambda
}(A_{2})\rightarrow H_{\lambda }(A_{1}\cup A_{2})\rightarrow H_{\lambda
-1}(A_{1}\cap A_{2})=0.$

The proof is by induction on $m$. In fact, we shall prove that Proposition A$%
_{(m,\lambda )}$ together with Proposition A$_{(m,\lambda +1)}$ implies
Proposition A$_{(m+1,\lambda ).}$

Suppose $F=\{A_{0},...,A_{m}\}$ is a finite collection of $m+1$ open subsets
of $X$ such that for $1\leq j\leq m+1$ and any subfamily $F^{^{\prime
}}\subset F$ of size $j,$ $\ \ H_{m-j+\lambda }(\bigcap F^{^{\prime }})=0.$
We will prove that $H_{m-1+\lambda }(\bigcup F)=0.$ Let us first prove,
using Proposition A$_{(m,\lambda +1)},$ that $H_{m-2+\lambda }(A_{1}\cup
...\cup A_{m})=0.$ This is so because for any subfamily $F^{^{\prime
}}\subset \{A_{1},...,A_{m}\}$ of size $j,$ $\ 1\leq j\leq m$, we have $\
H_{m-1-j+(\lambda +1)}(\bigcap F^{^{\prime }})=0.$

Let us consider the Mayer-Vietoris exact sequence of the pair $((A_{1}\cup
...\cup A_{m}),A_{0}):$

\begin{center}
$0=H_{m-1+\lambda }(A_{0})\oplus H_{m-1+\lambda }(A_{1}\cup ...\cup
A_{m})\rightarrow H_{m-1+\lambda }(\bigcup F)\rightarrow H_{m-2+\lambda
}(A_{0}\cap (A_{1}\cup ...\cup A_{m}))=0.$
\end{center}

Since by hypothesis, $H_{m-1+\lambda }(A_{0})=0,$ in order to conclude the
proof of Proposition A$_{(m+1,\lambda )}$ it is sufficient to prove that $%
H_{m-2+\lambda }(A_{0}\cap (A_{1}\cup ...\cup A_{m}))=0.$ For that purpose,
let $G=\{B_{1},...,B_{m}\}$ be the family of open subsets of $X$ given by $%
B_{i}=A_{0}\cap A_{i}$, $1\leq i\leq m.\ $Note that for any subfamily $%
G^{^{\prime }}\subset G$ of size $j,$ $\ 1\leq j\leq m$, $H_{m-1-j+\lambda
}(\bigcap G^{^{\prime }})=0.$ This is so because the homology group $\
H_{m-1-j+\lambda }(\bigcap G^{^{\prime }})=H_{m-(j+1)+\lambda }(\bigcap
F^{^{\prime }})=0,$ where $F^{^{\prime }}$ is the corresponding subfamily of 
$F$ of size $j+1$. Then by Proposition A$_{(m,\lambda )},$ $0=H_{m-2+\lambda
}(\bigcup G)=H_{m-2+\lambda }(A_{0}\cap (A_{1}\cup ...\cup A_{m})).$ \ This
completes the proof.\bigskip 
%TCIMACRO{\TeXButton{End Proof}{\qedsymbol}}%
%BeginExpansion
\qedsymbol%
%EndExpansion

We now give the Topological Berge's Theorem. See \cite{BM}.

\medskip

\noindent \textbf{Theorem B}$_{(m,\lambda )}.$\ \ \textit{Let }$%
F=\{A_{1},...,A_{m}\}$\textit{\ be a family of open subsets of a topological
space }$X$\textit{\ and let }$\lambda \geq 0$\textit{\ be an integer.
Suppose that}

\textit{\noindent a) }$H_{m-2+\lambda }(\bigcup F)=0;$

\textit{\noindent b) for }$1\leq j\leq m-1$\textit{\ and any subfamily }$%
F^{^{\prime }}\subset F$\textit{\ of size }$j,$

\begin{center}
$H_{m-2-j+\lambda }(\bigcap F^{^{\prime }})=0.$
\end{center}

\textit{\noindent Then \ }

\begin{center}
$H_{\lambda -1}(\bigcap F)=0.$
\end{center}

\noindent \textbf{Proof}. The proof is by induction. Theorem B$_{(2,0)}$
claims that two nonempty open sets whose union is connected must have a
point in common and Theorem B$_{(2,\lambda )}$ is just the statement of the
exactness of the Mayer-Vietoris sequence: \ $0=H_{\lambda }(A_{1}\cup
A_{2})\rightarrow H_{\lambda -1}(A_{1}\cap A_{2})\rightarrow H_{\lambda
-1}(A_{1})\oplus H_{\lambda -1}(A_{2})=0.$

Let us prove that Theorem B$_{(m,\lambda )}$ implies Theorem B$%
_{(m+1,\lambda ).}$ Let $F=\{A_{0,}A_{1},...,A_{m}\}$ as in Theorem B$%
_{(m+1,\lambda )}.$ That is, $H_{m-1+\lambda }(\bigcup F)=0$ and for $1\leq
j\leq m$\textit{\ \ }and any subfamily\textit{\ }$F^{^{\prime }}\subset F$%
\textit{\ }of size \textit{\ }$j$ we have $H_{m-1-j+\lambda }\mathit{(}%
\bigcap \mathit{F}^{^{\prime }}\mathit{)}=0$. Let $G=\{B_{1},...,B_{m}\},$
where $B_{i}=A_{0}\cap A_{i},$ $1\leq i\leq m.$ \ In order to prove that $\
H_{\lambda -1}(\bigcap F)=H_{\lambda -1}(\bigcap G)=0,$ it is enough to show
that the family $G=\{B_{1},...,B_{m}\}$ satisfies properties a) and b) of
Theorem B$_{(m,\lambda )}.$

Proof of a). \ We need to prove that $H_{m-2+\lambda }(A_{0}\cap (A_{1}\cup
...\cup A_{m-1}))=H_{m-2+\lambda }(\bigcup G)=0.$ Note that $H_{m-1+\lambda
}(\bigcup F)=0$ and $H_{m-2+\lambda }(A_{0})=0.$ Furthermore, by Proposition
A$_{(m,\lambda )},$ for the family $\{A_{1},...,A_{m}\},$ we have that $%
H_{m-2+\lambda }(A_{1}\cup ...\cup A_{m})=0.$ Thus the conclusion follows
from the Mayer-Vietoris exact sequence of the pair $(A_{0},(A_{1}\cup
...\cup A_{m}));$

\begin{center}
$0=$ $H_{m-1+\lambda }(A_{0}\cup ...\cup A_{m})\rightarrow H_{m-2+\lambda
}(A_{0}\cap (A_{1}\cup ...\cup A_{m}))\rightarrow H_{m-2+\lambda
}(A_{0})\oplus H_{m-2+\lambda }(A_{1}\cup ...\cup A_{m})=0.$
\end{center}

Proof of b). \ For $1\leq j\leq m-1$ and any subfamily $G^{^{\prime
}}\subset G$ of size $j,$ $\bigcap G^{^{\prime }}=\bigcap F^{^{\prime }},$
where $F^{^{\prime }}\subset F$ has size $j+1.$ Thus $H_{m-1-(j+1)+\lambda
}(\bigcap F^{^{\prime }})=H_{m-2-j+\lambda }(\bigcap G^{^{\prime }})=0.$
This completes the proof of Theorem B$_{(m+1,\lambda )}.$%
%TCIMACRO{\TeXButton{End Proof}{\qedsymbol}}%
%BeginExpansion
\qedsymbol%
%EndExpansion

\medskip

\bigskip

We now state our main theorem.

\medskip

\noindent \textbf{Topological Helly Theorem. }\textit{Let }$F$\textit{\ be a
finite family of open subsets of a topological space }$X$\textit{. Let }$d>0 
$\textit{\ be and integer such that }$H_{i}(U)=0$\textit{\ for }$i\geq d$%
\textit{\ \ and every open subset }$U$ \textit{of }$X.$\textit{\ }

\textit{Suppose that}

\begin{center}
$H_{d-j}(\bigcap F^{^{\prime }})=0$
\end{center}

\noindent \textit{for any subfamily }$F^{^{\prime }}\subset F$\textit{\ of
size }$j,$\textit{\ }$1\leq j\leq d+1.$ \ \textit{Then \ }

\begin{center}
$\bigcap F\neq \varnothing .$
\end{center}

\noindent \textit{Furthermore, }$\bigcap F$\textit{\ is acyclic.}

\noindent \textbf{Proof. \ }Suppose the size of $F$ is $m$. Take an integer $%
3\leq n\leq m+1.$ Using Theorem B$_{(n-1-\lambda ,\lambda )}$, from $\lambda
=0$ up to $\lambda =n-3,$ we can prove the following: \ 

\noindent Claim C$_{n}.$ Suppose that for every $1\leq j\leq n-1$ and any
subfamily $F^{^{\prime }}\subset F$ of size $j,$

\begin{center}
$H_{n-3}(\bigcup F^{^{\prime }})=0,$ and
\end{center}

\noindent for every $\ 1\leq j\leq n-2$ and any subfamily $F^{^{\prime
}}\subset F$ of size $j,$

\begin{center}
$H_{n-j-3}(\bigcap F^{^{\prime }})=0.$
\end{center}

\noindent Then for every $\ 1\leq j\leq n-1$ and any subfamily $F^{^{\prime
}}\subset F$ of size $j,$

\begin{center}
$H_{n-j-2}(\bigcap F^{^{\prime }})=0.$
\end{center}

Assume now $H_{\ast }(U)=0$ for every $\ast \geq d$ \ and every open $%
U\subset X$\ and suppose that $\ d\leq n-3.$ By repeating the use of Claim C$%
_{n},$ from $n=d+3$ up to $n=m+1,$ we obtain that $\bigcap F\neq \varnothing
.$

Arguing as above and using Theorem B$_{(m-1-\lambda ,\lambda )}$, from $%
\lambda =0$ up to $\lambda =m-3,$ we obtain that $H_{0}($ $\bigcap F)=0.$
The conclusion of acyclicity can be achieved by repeating the use of Theorem
B$_{(n,\lambda )},$ $2\leq n\leq m-1,$ $1\leq \lambda \leq m-3.$ Note that
in the case $m=d+2,$ our argument does not produce $H_{d-1}($ $\bigcap F)=0$%
, so we need to continue to $\lambda =m-2.$ This concludes the proof of our
main theorem.%
%TCIMACRO{\TeXButton{End Proof}{\qedsymbol}}%
%BeginExpansion
\qedsymbol%
%EndExpansion

\bigskip

For completeness, we include here a Topological Breen's Theorem.

\medskip

\noindent \textbf{Theorem} $\Sigma _{m}.$ \ \textit{Let }$%
F=\{A_{1},...,A_{m}\}$\textit{\ be a family of open subsets of a topological
space }$X$\textit{. Suppose that for }$1\leq j\leq m$\textit{\ \ and any
subfamily }$F^{^{\prime }}\subset F$\textit{\ of size }$j,$

\begin{center}
$H_{j-2}(\bigcup F^{^{\prime }})=0\mathit{.}$
\end{center}

\noindent \textit{Then}

\begin{center}
$\bigcap F\neq \varnothing \mathit{.}$
\end{center}

\noindent \textbf{Proof}. The proof is by induction. Theorem $\Sigma _{2}$
claims that two nonempty open sets whose union is connected must have a
point in common. Suppose Theorem $\Sigma _{m}$ is true and let $%
F=\{A_{0},A_{1},...,A_{m}\}$ be a family of open subsets of $X$ \ such that
for $1\leq j\leq m$ \ and any subfamily $F^{^{\prime }}\subset F$ of size $%
j, $ $\ H_{j-2}(\bigcup F^{^{\prime }})=0.$

Let us prove first that for any subfamily $F^{^{\prime }}\subset $ $%
\{A_{2},...,A_{m}\}$ of size $j,$ $0\leq j\leq m-1,$

\begin{center}
$H_{j-1}((A_{0}\cap A_{1})\cup \bigcup F^{^{\prime }})=0.$
\end{center}

To do so, simply consider the Mayer-Vietoris exact sequence of the pair $%
(A_{0}\cup \bigcup F^{^{\prime }},A_{1}\cup \bigcup F^{^{\prime }}):$

$0=H_{j}(A_{0}\cup A_{1}\cup \bigcup F^{^{\prime }})\rightarrow
H_{j-1}((A_{0}\cap A_{1})\cup \bigcup F^{^{\prime }})\rightarrow H_{j-1}$ $%
(A_{0}\cup \bigcup F^{^{\prime }})\oplus H_{j-1}(A_{1}\cup \bigcup
F^{^{\prime }})=0.$

This implies that the family $\{A_{0}\cap A_{1},A_{2},...,A_{m}\}$ satisfies
the hypothesis of Theorem $\Sigma _{m},$ and therefore by induction that $%
A_{0}\cap A_{1}\cap ...\cap A_{m}\neq \varnothing .$ This completes the
proof of this theorem.%
%TCIMACRO{\TeXButton{End Proof}{\qedsymbol}}%
%BeginExpansion
\qedsymbol%
%EndExpansion

\medskip

As an immediate consequence, we have the following theorem:

\medskip

\noindent \textbf{Topological Breen Theorem}.\textbf{\ }\textit{Let }$F$%
\textit{\ be a finite family of open subsets of a topological space }$X$%
\textit{. Let }$d>0$\textit{\ be and integer such that }$H_{\ast }(U)=0$%
\textit{\ \ for }$\ast \geq d$\textit{\ \ and every open subset }$U$ of $X$%
\textit{.}

\textit{Suppose that}

\begin{center}
$H_{j-2}(\bigcup F^{^{\prime }})=0$
\end{center}

\noindent \textit{for }$1\leq j\leq d+1$\textit{\ and} \textit{any subfamily 
}$F^{^{\prime }}\subset F$\textit{\ of size }$j.$ \textit{Then \ }

\begin{center}
$\bigcap F\neq \varnothing .$

\medskip
\end{center}

\noindent

\bigskip

\noindent \textbf{Remark}. The corresponding theorems for \v{C}ech
cohomology groups are also true. Furthermore, the theorems in this section
are true for a class of sets $A_{i}$ in which the Mayer-Vietoris exact
sequence of the pair $(A_{0},(A_{1}\cup ...\cup A_{m}))$ is exact. For
example, when $X$ is a polyhedron and every $A_{i}\subset X$ is a
subpolyhedron.

\bigskip \bigskip

\section{Transversal Theorems}

\subsection{Preliminary Lemmas}

We start with some notation.

Let $G(n,d),$\ be\textit{\ the Grassmannian space }of all\textit{\ }$n$%
\textit{-}planes in\textit{\ }$\mathbb{R}^{d}$\textit{\ }through the origin
and let $M(n,d)$ be the space of all affine $n$-planes in $\mathbb{R}^{d}$
as an open subset of $G(n+1,d+1).$

Let $F$ be a collection of nonempty convex sets in euclidean $d$-space $%
\mathbb{R}^{d}$ and let $0\leq n<d$ \ be an integer. \ We denote by $%
T_{n}(F)\subset $ $G(n,d)$ the topological space of all $n$-planes in $%
\mathbb{R}^{d}$ transversals to $F$; that is, the space of $n$-planes that
intersect all members of $F$. We say that $F$ is \textit{separated} if for
every $2\leq n\leq d$ and every subfamily $F^{^{\prime }}$ $\subset F$ of
size $n$, there is no $(n-2)$-plane transversal to $F^{^{\prime }}.$

Let $F=\{A^{1},...,A^{n}\}$ be a collection of closed subsets of a metric
space $X$ and let $\epsilon >0$ be a real number. We denote by $F_{\epsilon
}=\{A_{\epsilon }^{1},...,A_{\epsilon }^{n}\}$ the collection of open
subsets of $X$, where $A_{\epsilon }$ denotes the open $\epsilon $%
-neighborhood of $A\subset X.$

\medskip

\noindent \textbf{Lemma 3.1.1}. \textit{Let }$A$\textit{\ be a nonempty
convex set in }$\mathbb{R}^{d}$\textit{\ and let }$1\leq n<d.$\textit{\ Then 
}$T_{n}(\{A\})$\textit{\ is homotopically equivalent to }$G(n,d),$\textit{\
the Grassmannian space of all }$n$\textit{-planes in }$\mathbb{R}^{d}$%
\textit{\ through the origin.}

\noindent \textbf{Proof}. Let $\Upsilon :T_{n}(\{A\})\rightarrow G(n,d)$ be
given as follows:\ for every $H\in T_{n}(\{A\}),$ let $\Upsilon (H)$ be the
unique $n$-plane through the origin parallel to $H$. Then \ if $\ \Gamma \in
G(n,d),$ $\Upsilon ^{-1}(\Gamma )$ is homeomorphic to $\pi (A),$ where $\pi :%
\mathbb{R}^{d}\rightarrow \Gamma ^{\bot }$ is the orthogonal projection and $%
\Gamma ^{\bot }\in G(d-n,d)$ is orthogonal to $\Gamma .$ Since $\Upsilon $
has contractible fibers, it is a homotopy equivalence.%
%TCIMACRO{\TeXButton{End Proof}{\qedsymbol}}%
%BeginExpansion
\qedsymbol%
%EndExpansion

\medskip

\noindent \textbf{Lemma 3.1.2}. \textit{Let }$F=\{A_{1},A_{2},...,A_{n}\}$%
\textit{\ be a separated family of nonempty convex sets in }$\mathbb{R}^{d},$%
\textit{\ }$2\leq n\leq d,$\textit{\ and let }$n\leq m\leq d$\textit{\ \ be
an integer. Then }$T_{m-1}(F)$\textit{\ is homotopically equivalent to }$%
G(m-n,d-n+1).$

\noindent \textbf{Proof}. We start by proving that $T_{n-1}(F)$ is
contractible. For this purpose let $\Psi :A_{1}\times ...\times
A_{n}\rightarrow T_{n-1}(F)$ given by $\Psi ((a_{1},...,a_{n}))$ be equal to
the unique $(n-1)$-plane in $\mathbb{R}^{d}$ through $\{a_{1},...,a_{n}\},$
for every $(a_{1},...,a_{n})\in $ $A_{1}\times ...\times A_{n}.$ Note that $%
\Psi $ is well defined because $F$ is a separated family of sets.
Furthermore, if $H\in T_{n-1}(F),$ then $\Psi ^{-1}(H)=(H\cap A_{1})\times
...\times (H\cap A_{n})$ is contractible. This implies that $\Psi $ is a
homotopy equivalence and hence that $T_{n-1}(F)$ is contractible.

Let $E=\{(H,\Gamma )\mid H$ is a $(n-1)$-plane of $\mathbb{R}^{d},$ $\Gamma $
is a $(m-1)$-plane of $\mathbb{R}^{d}$ and $H\subset \Gamma \}.$ Then $%
\gamma :E\rightarrow M(n-1,d),$ given by the projection in the first
coordinate, is a classical fiber bundle with fiber $G(m-n,d-n+1)$. Now let $Y
$ $=\{(H,\Gamma )\in T_{n-1}(F)\times T_{m-1}(F)\ \mid H\subset \Gamma \}.$
Clearly, \ the restriction $\gamma \mid :Y\rightarrow T_{n-1}(F)$ is a fiber
bundle with fiber $G(m-n,d-n+1)$ and contractible base space $T_{n-1}(F).$
Therefore $\gamma \mid :Y\rightarrow T_{n-1}(F)$ is a trivial fiber bundle
and hence $Y$ is homotopically equivalent to $G(m-n,d-n+1).$

Consider now the projection $\pi :Y\rightarrow T_{m-1}(F).$ Note that for
every $\Gamma \in T_{m-1}(F),$ the fiber $\pi ^{-1}(\Gamma )$ is equal to $%
T_{n-1}(\{A_{1}\cap \Gamma ,...,A_{n}\cap \Gamma \}).$ By the first part of
this proof, the fibers of $\pi $ are contractible, hence $\pi $ is a
homotopy equivalence and $T_{m-1}(F)$ is homotopically equivalent to $%
G(m-n,d-n+1).$%
%TCIMACRO{\TeXButton{End Proof}{\qedsymbol}}%
%BeginExpansion
\qedsymbol%
%EndExpansion

\medskip \bigskip

\noindent \textbf{Lemma 3.1.3}. \textit{Let }$A,B$\textit{, and }$C$\textit{%
\ be three nonempty convex sets in }$\mathbb{R}^{d}$\textit{\ such that }$%
A\cap B=\varnothing .$\textit{\ Then}

\begin{center}
$H_{1}(T_{1}(\{A,B,C\}))=0.$\textit{\ \ \ }
\end{center}

\noindent \textbf{Proof}. Since $A\cap B\cap C=\varnothing ,$ by Theorem 3
of \cite{BMO}, $T_{1}(\{A,B,C\})$ has the homotopy type of the space $%
\mathcal{C}_{1}(\{A,B,C\})\subset \mathcal{C}_{2}^{1}$ \ of all affine
configurations of three points in the line, achieved by transversal lines to 
$\{A,B,C\}.$ Note now that the space of affine configuration of three points
in a line, $\mathcal{C}_{2}^{1}$, is $\mathbb{S}^{1}$ and note further that
since $A\cap B=\varnothing ,$ the space of all affine configurations of
three points in the line achieved by transversal lines to $\{A,B,C\}$ is a
subset $\mathbb{S}^{1}-\{\infty \},$ where $\infty \in \mathbb{S}^{1}$ is
the affine configuration in which the first and the second points coincide.
This implies that $H_{1}(T_{1}(\{A,B,C\}))=0.$\ 
%TCIMACRO{\TeXButton{End Proof}{\qedsymbol}}%
%BeginExpansion
\qedsymbol%
%EndExpansion

\medskip \bigskip

\noindent \textbf{Lemma 3.1.4}. \ \textit{Let }$F=\{A^{1},...,A^{d+1}\}$ 
\textit{be a separated family of closed convex sets in }$\mathbb{R}^{d}$. 
\textit{Suppose that } $H_{0}(T_{d-1}(F))=0.$ Then there is $\epsilon _{0}>0$
with the property that if \textit{\ }$0<\epsilon <\epsilon _{0},$ then%
\textit{\ } $H_{0}(T_{d-1}(F_{\epsilon }))=0.$

\noindent \textbf{Proof.} By Theorem 1 of \cite{BMO}, the space of
transversals $T_{d-1}(F)$ of a separated family of convex sets in $\mathbb{R}%
^{d}$ has finitely many components and each one of them is contractible. In
fact, each component corresponds precisely to a possible order type, of $d-1$
points in affine $(d-1)$-space, achieved by the transversal hyperplanes when
they intersect the family $F$. In our case, since $H_{0}(T_{d-1}(F))=0,$ we
have that $T_{d-1}(F)$ is contractible and that the transversal hyperplanes
intersect the family $F$ consistently with a precise order type $\Omega .$

Suppose now the lemma is not true, then there exist an order type $\Omega
_{0}$, different from $\Omega ,$ and a collection of hyperplanes $H_{i}$
that intersect $F_{\epsilon _{i}}$ consistently with the order type $\Omega
_{0.}$ Since we may assume that $\{\epsilon _{i}\}\rightarrow 0$ and $%
\{H_{i}\}\rightarrow H,$ where $H$ is a transversal hyperplane to $F$
consistently with the order type $\Omega _{0},$ we have a contradiction. 
%TCIMACRO{\TeXButton{End Proof}{\qedsymbol}}%
%BeginExpansion
\qedsymbol%
%EndExpansion

\medskip \bigskip

\subsection{Transversal Lines in the Plane}

A family of sets is called \textit{semipairwise disjoint} if, given any
three elements of $F$, two of them are disjoint.

\medskip

\noindent \textbf{Theorem 3.2.1. }\textit{Let }$F$\textit{\ be a
semipairwise disjoint family of at least 6 open \textit{c}onvex sets in }$%
\mathbb{R}^{2}.$\textit{\ Suppose that for every subfamily }$F^{^{\prime }}$%
\textit{\ }$\subset F$\textit{\ \ of \ size 5, }$T_{1}(F^{^{\prime }})\neq
\phi $\textit{\ and for every subfamily }$F^{^{\prime }}$\textit{\ }$\subset
F$\textit{\ \ of \ size 4, }$T_{1}(F^{^{\prime }})$\textit{\ is connected. \
Then }$T_{1}(F)\neq \varnothing $\textit{.}

\noindent \textbf{Proof}. Let $X$ be the space of all lines in $\mathbb{R}%
^{2}.$ Hence $H_{\ast }(U)=0$ for $\ast \geq 2$ \ and every open subset $%
U\subset X$. We are interested in applied the Topological Helly Theorem when 
$d=4$.\ Note first that $H_{3}(T_{1}(\{A\})=0$ for every $A\in F,$ \ and $%
H_{2}(T_{1}(\{A,B\})=0$ \ for $A\neq B\in F.$ By Lemma 3.1.3 and the fact
that $F$ is semipairwise disjoint, we have that $H_{1}(T_{1}(F^{^{\prime
}}))=0,$ for any subfamily $F^{^{\prime }}\subset F$ of size $3$. By
hypothesis, $H_{0}(T_{1}(\{F^{^{\prime }}\})=0,$ for any subfamily $%
F^{^{\prime }}\subset F$ of size $4,$ and $H_{-1}(T_{1}(\{F^{^{\prime
}}\})=0,$ for any subfamily $F^{^{\prime }}\subset F$ of size $5$. This
implies, by the Topological Helly Theorem, that $T_{1}(F)$ is nonempty.%
%TCIMACRO{\TeXButton{End Proof}{\qedsymbol}}%
%BeginExpansion
\qedsymbol%
%EndExpansion

\medskip \medskip \medskip

\subsection{Transversal Lines in $3$-Space}

In this section we study transversal lines to families of convex sets in $%
\mathbb{R}^{3}.$

\medskip

\noindent \textbf{Theorem 3.3.1.} \textbf{\ }\textit{Let }$F$\textit{\ \ be
a pairwise disjoint family of at least }$6$\textit{\ open, convex sets in }$%
\mathbb{R}^{3}$\textit{.\ Suppose that for any subfamily }$F^{^{\prime }}$%
\textit{\ }$\subset F$\textit{\ \ of \ size }$5$\textit{, }$%
T_{1}(F^{^{\prime }})\neq \varnothing ,$\textit{\ and for any subfamily }$%
F^{^{\prime }}$\textit{\ }$\subset F$\textit{\ \ of \ size 4, }$%
T_{1}(F^{^{\prime }})$\textit{\ is connected. \ Then, }$T_{1}(F)\neq
\varnothing $\textit{.}

\noindent \textbf{Proof}. Let $X$ be the space of all lines in $\mathbb{R}%
^{3},$ hence $X$ is an open $4$-dimensional manifold and therefore $H_{\ast
}(U)=0$ for $\ast \geq 4$ \ and every open subset $U\subset X$. We are
interested in applied the Topological Helly Theorem for $d=4$.\ By Lemma
3.1.1, $H_{3}(T_{1}(\{A\})=0,$ for every $A\in F,$ since $T_{1}(\{A\})$ has
the homotopy type of $G(1,3)=\mathbb{RP}^{2}.$ By Lemma 3.1.2, $\
H_{2}(T_{1}(\{A,B\}))=0,$ for every $A\neq B\in F.$ By Lemma 3.1.3 and the
fact that $F$ is pairwise disjoint, we have that $H_{1}(T_{1}(F^{^{\prime
}}))=0,$ for any subfamily $F^{^{\prime }}\subset F$ of size $3$. By
hypothesis, $H_{0}(T_{1}(\{F^{^{\prime }}\})=0,$ for any subfamily $%
F^{^{\prime }}\subset F$ of size $4,$ and $H_{-1}(T_{1}(\{F^{^{\prime
}}\})=0,$ for any subfamily $F^{^{\prime }}\subset F$ of size $5$. This
implies, by the Topological Helly Theorem, that $T_{1}(F)$ is nonempty. 
%TCIMACRO{\TeXButton{End Proof}{\qedsymbol}}%
%BeginExpansion
\qedsymbol%
%EndExpansion

\medskip

\noindent

\subsection{Transversal Hyperplanes}

This section is devoted to stating and proving a theorem concerning
transversal hyperplanes to families of separated convex sets in d-space.

\medskip

\noindent \textbf{Theorem 3.4.1}. \textit{Let }$F$\textit{\ be a separated
family of at least }$d+3$ \textit{closed, convex sets in }$\mathbb{R}^{d}.$ 
\textit{Suppose that for any subfamily }$F^{^{\prime }}$\textit{\ }$\subset
F $\textit{\ \ of \ size }$d+2$\textit{, }$T_{d-1}(F^{^{\prime }})\neq
\varnothing $\textit{\ and for any subfamily }$F^{^{\prime }}$\textit{\ }$%
\subset F$\textit{\ \ of \ size }$d+1$\textit{, }$T_{d-1}(F^{^{\prime }})$%
\textit{\ is connected. \ Then }$T_{d-1}(F)\neq \varnothing $\textit{.}

\noindent \textbf{Proof}. Let us first prove the theorem for a separated
family of open convex sets. We are going to use the Topological Helly
Theorem. Let $X$ be the space of all hyperplanes of $\mathbb{R}^{d}.$ Note
that $H_{\ast }(U)=0$ for $\ast \geq d$ \ and every open subset $U\subset X.$
In particular, $H_{\ast }(U)=0$ for every $\ast \geq d$ $+1.$

By Lemma 3.1.2, for every subfamily $F^{^{\prime }}$ $\subset F$ \ of \ size 
$j$, $1\leq j\leq d,$ $T_{d-1}(F^{^{\prime }})$ is homotopically equivalent
to $G(d-j,d-j+1)$ and hence $H_{d-j+1}(T_{d-1}(F^{^{\prime }}))$ $\ =$ $%
H_{d-j+1}(\bigcap \{T_{d-1}(\{A\})\mid A\in F^{^{\prime }}\})=0.$
Furthermore, by hypothesis, the same is true for $j=d+1$ and $j=d+2.$
Consequently, by our Topological Helly Theorem, $T_{d-1}(F)\neq \varnothing $%
.

By Lemma 3.1.4, there is $\epsilon >0,$ such that $F_{\epsilon }$\textit{\ }%
is a separated family of open, convex sets in\textit{\ }$\mathbb{R}^{d}$ and
for any subfamily $F_{\epsilon }^{^{\prime }}\mathit{\ }\subset F_{\epsilon }%
\mathit{\ \ }$of size $d+1,_{\text{ \ }}T_{d-1}(F_{\epsilon }^{^{\prime }})$
is connected. By the above, this implies that $T_{d-1}(F_{\epsilon }))\neq
\varnothing .$ Hence, by completeness of the Grassmannian spaces, $%
T_{d-1}(F)\neq \varnothing $.%
%TCIMACRO{\TeXButton{End Proof}{\qedsymbol}}%
%BeginExpansion
\qedsymbol%
%EndExpansion

\end{document}